\documentclass[draft, leqno, oneeqnum, onethmnum, translated]{tvp}

\def\FnMarkWithoutParanthesis{
\def\@makefnmark{\hbox
to 7pt{$^{{\@thefnmark}}$}} 
\long\def\@makefntext##1{\parindent 1em\noindent
\hbox to 2.2em{\hss$^{{\@thefnmark}}$}\ ##1}
}

\def\FnMarkWithParanthesis{
\gdef\@makefnmark{\hbox
to 7.8pt{$^{{\@thefnmark})}$}}  
\long\def\@makefntext##1{\parindent 1em\noindent
\hbox to 2.2em{\hss$^{{\@thefnmark})}$}\ ##1}
}

\newcommand{\otk}{``\ignorespaces}

\newcommand{\zak}{''}

\newcommand{\llla}[1]{\mathop{\,\hbox to #1{\rightarrowfill}\,}}

\renewcommand{\le}{\leqq}
\renewcommand{\ge}{\geqq}

\newcommand{\ConjectureName}{Conjecture}

\newtheorem{def-prop}{Definition-proposition}
\newtheorem{prop-def}{Proposition-definition}

\newtheorem{def-propi}{Definition-proposition}

\newtheorem{def-propsi}{Definition-proposition}[section]

\newtheorem{notationnon}{Notation}

\newcommand{\Translated}{Translated from Russian Journal}

\font\sixrm=cmr6

\newcounter{VolumeNumber}
\newcounter{IssueNumber}

\newcommand{\PaperHead}{{\scriptsize\hbox to\hsize{{{\sixrm T}{\tiny\sc HEORY}
{\sixrm P}{\tiny\sc ROBAB}.\ {\sixrm A}{\tiny\sc PPL}.}\hfill{\sixrm\Translated}}
\hbox to\hsize{\sixrm{Vol.\ \theVolumeNumber, No.\ \theIssueNumber}\hfill}}\vspace{9pt}}

\newcommand{\remno}[1]{{\em Remark\/ }#1. }

\newcommand{\exampno}[1]{{\em Example\/ }#1. }

\newcommand{\kb}{{\upshape, }}
\newcommand{\kbt}{{\upshape; }}

\newcommand{\1}{\hbox{$1$\kern-3pt{\rm I}}\kern.5pt } 
\newcommand{\0}{\hbox{$0$\kern-2.8pt\raisebox{-.26pt}\hbox{\vbox%
               {\hrule height5.8pt width0.4pt}}}\hskip3.4pt} 

\newlength{\smallcolsep}
\newlength{\largecolsep}
\newcommand{\cl}{\colon\,}

\newcommand{\beq}{\begin{equation}}
\newcommand{\eeq}{\end{equation}}
\newcommand{\bea}{\begin{eqnarray}}
\newcommand{\eea}{\end{eqnarray}}
\newcommand{\bnz}{\begin{eqnarray*}}
\newcommand{\enz}{\end{eqnarray*}}

\newcommand{\bE}{\mathbf{E}}

\newcommand{\bP}{\mathbf{P}}
\newcommand{\bR}{\mathbf{R}}

\newcommand{\ccF}{\mathcal{F}}

\newcommand{\sH}{\mathbb{H}}

\newcommand{\rrd}{{\rm d}}

\usepackage{amssymb}
\usepackage{amsmath}
\usepackage{exscale}
\def\E{{\bf E}}

\author{Ulyanov~V.\,V.\thanks{Lomonosov Moscow State University, Faculty of Computational Mathematics and Cybernetics; National Research University Higher School of Economics (HSE), Moscow, Russia; e-mail: vulyanov@cs.msu.ru}}

\title{ON PROPERTIES OF POLYNOMIALS IN RANDOM ELEMENTS\thanks{This research was supported by RSCF
  14-11-00196.}}

\begin{document}


\maketitle



\begin{abstract} The paper deals with different properties of polynomials in random elements: bounds for characteristics functionals of polynomials, stochastic generalization of the Vinogradov mean value theorem, characterization problem, bounds for probabilities to hit the balls. These results cover the cases when the random elements take values in finite as well as infinite dimensional Hilbert spaces.
\end{abstract}


\pagestyle{myheadings} 
\thispagestyle{plain}  

\setcounter{page}{1}



\section{Introduction}\ Let $X$~be a random element with
values in a real separable Hilbert space ${\sH}$. We consider a sequence  $X_1, X_2, \dots$ of
independent copies of $X$. Put
$$
S_n = n^{-1/2}(X_1+\cdots +X_n), \quad g_n(t) = \E \exp\{itf(S_n)\},
$$
where $f$~is a measurable map from ${\sH}$ to the real line ${\bR}$.

In order to find estimates for the rate of convergence in functional limit theorems and
to construct asymptotic expansions, we apply inequalities for
  $|g_n(t)|$ of the type
\begin{equation}
|g_n(t)| \le c(n^{-A} + |t|^{-D}),
\label{u1}
\end{equation}
where $c, A$ и  $D$~are some constants. One can derive estimate
  (\ref{u1}) using the so-called symmetrization inequality (see
   .~[1] and~[2], in particular, Corollary~3.26 and the remark after
it in~[2] concerning the paper of Weyl~[3]; see\ also~[4]). For example, when ${\sH} = \bR$ this approach leads to the following inequality (cf.\ [5, Lemmas~3.1, 4.1, 4.2]). 

Set
\begin{equation}
f(x) = x^m + \alpha_{m-1}x^{m-1} + \cdots +\alpha_0, \qquad m\ge 2.
\label{2}
\end{equation}

Then we have (\ref{u1})   with $D = (m\cdot 2^m)^{-1}$ and an arbitrary $A > 0$ for~$t$, satisfying the condition
\begin{equation}
|t|\le c_1n^{(m-1)/2-\varepsilon},
\label{3}
\end{equation}
where $\varepsilon > 0$.

It is well known that if a random vector $X$ has a discrete distribution, then $g_n(t,a)$
is an almost periodic function. Therefore,
$$
\limsup_{t\rightarrow\infty}|g_n(t,a)| = 1.
$$
Thus, we have (\ref{u1})   provided
  $t$ is restricted to a finite interval depending on~$n$ and~$f$.

In Section~2 we give the bounds for characteristic functions $g_n(t,a)$, which improve~(\ref{u1}) so that
  ~$D$ increases essentially and becomes proportional to $1/m$. It follows from Theorem~4 with two-sided bounds for  $g_n(t,a)$, that this order of~$D$ is right one. The improvements of~(\ref{u1}) in Section~2 are proved provided that some conditions are met concerning the type of distribution of~$X$.
  A desire to remove these conditions leads to the necessity of stochastic generalization of the famous Vinogradov mean value theorem. The stochastic generalization see in Section~3. In particular, if $X$ takes values
   $1,2,\dots, P$,
   with equal probabilities, our estimate gives the same order with respect to~$P$ as Vinogradov's original result.
    The properties of quadratic forms are considered in Sections~4 and~5.
     They attract the attention of researchers in recent years.
      In Section~4 we give sufficient conditions in order a distribution of a finite quadratic form in independent identically distributed symmetric random variables define a distribution of the basic random variable uniquely. The stability theorem for quadratic forms is proved as well. In Section~5 the two-sided bounds are found for the density $p(u,a)$ of  $|Y-a|^2$, where $Y$ is an~${\sH}$-valued Gaussian random element. The bounds are precise in the sense that the ratio of upper bound to the lower one equals 8 for all sufficiently large values of $u$. Thus, the ratio does not depend on the parameters of the distribution of
 $|Y-a|^2$. These bounds imply two-sided bounds for the probabilities $\bP(|Y-a| > r).$

In this paper we use and discuss Yu.V.Prokhorov's results obtained jointly with coauthors in
 1995--2000 \ (see~[6]--[11]). Asymptotic behavior of quadratic and almost quadratic forms that appeared in mathematical statistics is considered in~[12]. See~[13] as well.

\section{Bounds for characteristic functions of polynomials in random elements}\
As noted in the introduction,~(\ref{u1}) could be improved so that $D$ becomes larger (see~[6],~[7]). The
improvement is obtained in the following cases:

a)~the distribution of $X$ has a non-degenerate discrete component;

b)~for some $n_0$ the distribution of $S_{n_0}$ has an absolutely continuous component;

c)~the characteristic function of $X$ satisfies some conditions \otk in average\zak, in particular, the conditions are met for a class of singular distributions.

Recall that if $\bP, \bP_1, \bP_2$ denote three probability measures in   ${\sH}$ and $\gamma$ is a positive
number, $0 < \gamma \le 1$, then $\bP_1$ is called  {\it a component} of  $\bP$ of {weight}~$\gamma$, provided
$$
{\bP} = \gamma\, {\bP_1} + (1-\gamma)\, {\bP_2}.
$$
In case a) the following theorem is true (see [6, Theorem 2]).

\begin{theorem} \label{th1.1}
Let $X_1, X_2, \dots$~be independent identically distributed $($i.i.d.$)$ random variables in~$\bR$ with a distribution~$F,$ admitting
a non-degenerate discrete component. Then for any $\varepsilon > 0$ and integer $m \ge 2$ there exist absolute constants
  $c_3, c_4$ and a constant $c_2, $  depending on $F,$ $\varepsilon$ and $m,$ such that\kb for all $n \ge 1$ and any  $t,$ in the range~{\rm(\ref{3}),} one has
\begin{equation}
\sup_{a\in\bR} |\bE\,\exp\{itf(S_n + a)\}|\le c_2 \bigg(n^{-c_3/(m^2 \ln m)} + \bigg(\frac{\sqrt{n}}{|t|}\bigg)^{c_4/m}\bigg),
\end{equation}
where $f$~is the polynomial defined in~{\rm(\ref{2})}.
\end{theorem}

Now we consider the case of continuous components. We assume that a random
vector
  $Y$ in ${\bR}^k$ has a distribution  ${\bP}_Y$ which admits an absolutely continuous component
of weight $\lambda$. Then ${\bP}_Y$ has an absolutely continuous component of weight $\gamma$, $0 < \gamma < \lambda$, such that the Lebesgue density of this component is bounded. Take $\gamma = \lambda/2$. Then (see~[14, p.~4] or [15, \S\,16]) this implies that
${\bP}_Y \ast {\bP}_Y$ has a component
of weight $(\lambda/2)^2$ with uniformly continuous density function on $\bR^k$, say, $r(x)$. If $r(a_0) = r_0 > 0$ for some $a_0$,
then there exists $\eta$ such that for all $x$, with $|x - a_0|\le \eta$, we have $r(x)\ge r_0/2$. Thus, ${\bP}_Y \ast {\bP}_Y$  admits a component which is uniformly
distributed on some ball.

In case b) the following theorem is true (see~[6, Theorem~3]).

\begin{theorem} \label{th1.2}
   Set for a vector $x = (x^{(1)}, \dots , x^{(k)})$
$$
f_0(x) = (x^{(1)})^m +\cdots + (x^{(k)})^m, \qquad m\ge 2.
$$
Let
  $X_1, X_2, \dots$ denote a sequence of i.i.d.\ random vectors in
  $\bR^k$. Suppose that for some
    $n_0$ the distribution of $X_1 + \cdots + X_{n_0}$   has an absolutely continuous
component. Let $\eta$ denote a real number such that the distribution of $X_1 + \cdots + X_{2n_0}$ has a component of weight $\delta,$ uniformly distributed on the ball $\{x\cl |x - a_0|\le \eta\}$. Then there exist the constants   $c_5, c_6,$ depending on $k, n_0, \delta, \eta$ and such that for all $n \ge 6\,n_0/\delta$ and $t\neq 0$ we have
\begin{eqnarray*}
\sup_{a\in \bR^k} |\bE\,\exp\{it\,f_0(S_n + a)\}|
 \le c_5(\exp \{-c_6n\} + |t|^{-k/m}).
\end{eqnarray*}
\end{theorem}

The expressions for   $c_5, c_6$ see\ in~[6, Theorem~3].
For components with differentiable densities Theorem~\ref{th1.2} could be refined (see [6, Theorem 4]).

We now consider the general polynomial
\begin{equation}
f(x) = \sum_{m_1=0}^m \cdots \sum_{m_k=0}^m \alpha(m_1, \dots , m_k)\,x_1^{m_1}\cdots x_k^{m_k},
\label{4}
\end{equation}
where $\alpha(0, \dots , 0) = 0$. Let $M$ stand for the {\it degree} of $f(x)$, that is,
$$
M = \max \{m_1 + \cdots + m_k\cl \alpha(m_1, \dots , m_k) \neq 0\}.
$$
Clearly, $km \ge M$. Put
$$
\alpha_{*} = \max |\alpha(m_1, \dots , m_k)| > 0,
$$
where the maximum is taken over all $m_1, \dots , m_k$ such that
$$
0 \le m_1, \dots , m_k \le m, \quad m_1 + \cdots + m_k = M.
$$
The following theorem is true (see\ [6, Theorem~5]).

\begin{theorem} \label{th1.3}
Let $Z$~be a random vector in ${\bR}^k$ with independent standard normal
coordinates $Z^{(i)}$\kb $i = 1, \dots , k$. If $f(x)$\kb is defined by~{\rm(\ref{4}),} with $M \ge 2,$ then
$$
\sup_{a\in \bR^k} |\bE\,\exp\{itf(Z + a)\}|
 \le c_7|\alpha_{*}t|^{-1/m}\ln^{c_8}(2 + |\alpha_{*}t|)
$$
for some constants $c_7, c_8,$ depending on $k, M, m$.
\end{theorem}

The expressions for $c_7, c_8$ see\ in [6, Theorem~5].

Notice that until now we gave the
  {\it upper} bounds for characteristic functionals. In some cases it is possible to prove two-sided bounds. Moreover, the orders of the upper and lower bounds with respect to~$t$ coincide. Set
$$
f_1(x) = x^{(1)} \cdots x^{(k)}, \qquad k\ge 2.
$$
Let $Z$, as above, be a random vector in ${\bR}^k$ with independent standard normal
coordinates.
Applying Theorem~3, we obtain
$$
|\bE\,\exp\{it f_1(Z)\}| \le c(k)|t|^{-1}\ln^{k-1}(2 + |t|)
$$
only. In fact, one can get not only better upper bound but prove a lower bound of the same order in
 ~$t$ as well (see\ [6, Theorem~6]).

\begin{theorem} \label{th1.4}
For any $k \ge 2$ and all $|t|\ge 1$ we have
$$
l_k|t|^{-1}\ln^{k-2}|t| \le |\bE\,\exp\{it f_1(Z)\}|  \le L_k |t|^{-1}\ln^{k-2}(2 + |t|)
$$
with some $l_k, L_k,$ depending on~$k$ only.
\end{theorem}

The expressions for $l_k, L_k$ see\ in [6, Theorem~6].

Now we consider the case when the characteristic function of
 ~$X$ satisfies some \otk averaged\zak \ conditions. Let $g_X(t) = \bE\,\exp \{itX\}$ be the characteristic function of~$X$. Put for $T > 0$
$$
\phi_X(T) = \int_{-T}^T|g_X(t)|\,dt.
$$
The following theorem is true (see\ [7, Theorem~1]).

\begin{theorem} \label{th1.5}
Let $X$~be a random variable such that $|X|\le 1$. Assume that there
exists a nondecreasing positive function $\phi(t)$ on $(0, \infty)$ and a constant $\varepsilon:\, 0\le \varepsilon < 1/m$\kb such that for all $t \ge 1$ and $b \ge 1$
\begin{equation}
\phi_X(bt) \le b^{\varepsilon} \phi(t).
\label{5}
\end{equation}
Then for  $|t| \ge 1$ we have
\begin{equation}
\Big|\bE\, \exp \{itf(X)\}\Big| \le c\phi(|t|)|t|^{-1/m},
\label{6}
\end{equation}
where   $f(x)$ is defined by {\rm(\ref{2})} and the constant $c$ does not depend on~$t$.
\end{theorem}

The expressions for $c$ see\ [7, relations~(2.5) and~(2.10)].

\remno{1} If the characteristic function $g_X(t)$ satisfies
\begin{equation}
\int_{-\infty}^{\infty}|g_X(t)|\,dt = A < \infty,
\label{7}
\end{equation}
we may choose, in Theorem~\ref{th1.5},   $\varepsilon = 0$ and $\phi(t) = A$. Thus, assuming condition~(\ref{7}) we get
\begin{equation}
\Big|\bE\,\exp \{itf(X)\}\Big| \le c A |t|^{-1/m}.
\label{8}
\end{equation}
Moreover, if (\ref{7}) holds we drop the assumption $|X|\le 1$ (see\ [6, Lemma~7 and bounds~(22),~(23)]).

\remno{2} One can show (see\ [7, Corollary~3]), the distribution function of~$X$,
which
satisfies the conditions of Theorem~\ref{th1.5} is necessarily a Lipschitz function of order
  $1 - \varepsilon$.

Condition (\ref{5}) is satisfied by many absolutely continuous distributions (cf.\ [16, Ch.~VI, \S\,3 ])
as well as by some singular distributions
  (see~[17]--[19] and [16, Ch.~XII, \S\S\,10, 11]).

We give an example of singular distribution of~$X$, such that its convolutions of any order are
singular as well. However, the behavior of the characteristic function of polynomial mapping
  $f(X)$ is very similar to the case when~$X$ is absolutely
continuous (see\ Theorem~\ref{th1.2}) or even Gaussian (see\ Theorem~\ref{th1.3}).

We consider the Cantor distribution with characteristic function
\begin{equation}
L(t) = \prod_{j=1}^{\infty} \cos (2\pi\cdot 3^{-j}t).
\label{9}
\end{equation}
The Cantor distribution satisfies the Cram\`{e}r condition (see [17]), i.e.,
\begin{equation}
\limsup_{|t|\rightarrow \infty} |L(t)| < 1.
\label{10}
\end{equation}
In [7, Theorem~2] it was proved  a sharper bound than~(\ref{10}), namely, for all
$|t| \ge 8.5$ we have
$$
|L(t)| \le e^{-0.027}.
$$
Moreover, in fact, not only is the Cantor distribution singular but its $k$-fold convolutions of any order $k$ are
singular as well (see [7, Theorem~2]). The convolutions have characteristic functions $L^k(t)$. At the same time the behavior of the characteristic function of the Cantor distribution after polynomial mapping is similar to the case of absolutely continuous distributions.

The following theorem is true (see\ [7, Corollary~5]).

\begin{theorem} \label{th1.6} Let $f(x)$~be a polynomial defined in~{\rm(\ref{2})} with $m\ge 2$, and let  $\varepsilon : \,0< \varepsilon < 1/m$. Let~$X$
denote a random variable with singular distribution and characteristic function
  $L^k(t)$\kb
where $k$~is any integer such that $k\ge 0.027^{-1}\ln (2/(3^{\varepsilon} - 1))$. Then
$$
\bE\,\exp \{itf(X)\} = O (|t|^{\varepsilon - 1/m}) \quad \mbox{as}\ |t| \rightarrow \infty.
$$
\end{theorem}

The proofs of Theorems~1--4
are based on inequalities for trigonometric sums and integrals. If the distribution of $X$ has a non-degenerate discrete component we use estimates of the Vinogradov-type for trigonometric sums (see\ [20, Statement~4.2 in \S\,4, Ch.~VIII]). If for some $n_0$ the distribution of $S_{n_0}$ has an absolutely continuous component, we apply the Vinogradov inequality for trigonometric integrals (see\ [21, Lemma~4, Ch.~2])
and its generalizations to the multidimensional case (see [22, Theorem~5 in \S\,3, Ch.~1]).
Our improvements
of estimate~(\ref{u1}) obtained in~[6]  are similar in form with the improvements which gave the Vinogradov
method for trigonometric sums comparing the Weyl results (see Introduction in [21]).

\section{Stochastic generalization of the Vinogradov mean value theorem}\
A desire to remove additional conditions on the type of the distribution X leads to the
necessity of stochastic generalization of the famous Vinogradov theorem on the mean
Theorems~1--4 are proved provided additional conditions on the type of the distribution~$X$ are met.
A desire to remove the conditions leads to the
necessity of stochastic generalization of the famous Vinogradov mean value theorem
  (see.~[21]). The stochastic generalization was obtained in~[8]. Introduce necessary notation.

Let $J_k(P)$ denote the number of the simultaneous diophantine equations
$$
\sum^k_{i=1} (x_i^j-y_i^j)=0,\qquad j=1,\dots,m,
$$
where $1\le x_i,$ $y_i\le P$ for $i=1,\dots, k.$

Since for any integer $x$ we have
$$
\int^1_0 \exp\{2\pi i \alpha x \}\,d\alpha=
\left\{\begin{array}{ll}
1, \quad & \mbox{if}\ x=0, \\
       0, \quad & \mbox{if}\ x \ne 0,
\end{array}\right.
$$
we get
$$
J_k(P)=\int^1_0\dots \int^1_0\bigg|\sum^P_{x=1}\exp\{2\pi if (x)\}
     \bigg|^{2k} \,d\alpha_m\cdots d\alpha_1,
$$
where $f(x) = \alpha_mx^m+\cdots+\alpha_1x.$

Estimates of the magnitude order of the growth of
   $  J_k(P)$
with increasing $P$ have important applications in analytic number theory. In 1934 I.M.Vinogradov proved the following theorem for
 $J_k(P)$.

\begin{theorem} Let $\tau > 0,$ $m>2 $~be integers with  $k\ge m\tau,$ and let $P\ge 1$. Then
$$
J_k(P)\le c_{\tau}\cdot P^{2k-\triangle(\tau)},
$$
where
$$
\triangle(\tau)=0.5 m(m+1)(1-(1-m^{-1})^\tau),\quad
c_{\tau}=(m\tau)^{6m\tau} (2m)^{4m(m+1)\tau}.
$$
\end{theorem}

The bounds for $ J_k(P)$ provide information about the size of the trigonometric Weyl sums
$$
F = \sum^P_{x=1} \exp\{2\pi if(x)\}
$$
and trigonometric sums involving functions to which there are reasonable polynomial approximations. Since there are numerous applications for such bounds Vinogradov was able to use his method with great success in many different problems of number theory. For later refinements of the Vinogradov mean value theorem see, e.g.,~[23]--[25], and the discussion and references therein.

It is obvious that if a random variable
 $S$ takes values   $1,2,\dots,P$ with the equal probability
  $1/P$, then
$$
\bE\,\exp\{2\pi if (S)\}=\frac{F}{P}.
$$
In Section~2 we showed that using the bounds for   $|F|$, known in the number theory, one can improve the upper bounds for
$|\bE\,\exp\{2\pi if (S_n)\}|,$ where $S_n$ is a normalized sum of i.i.d. random variables   $X_1,\dots,X_n$ such that the distribution of~$X_1$ has non-degenerate component
   (see~[6]).

Stochastic generalization of the Vinogradov mean value theorem which is of  independent interest, could be applied to find the upper bounds for
$|\bE\,\exp\{2\pi if (S_n)\}|$ without any additional conditions for the type of distribution of~$X_1$.

The following theorem is true (see
 \ [8, Theorem~2]).

\begin{theorem} Let $S$ be a random variable and $P\ge1,$ $m>2,$ $\tau\ge 1$ and $k$ denote the positive integers.  Set $f(x)=\alpha_m x^m+\cdots+\alpha_1 x$ and
\begin{eqnarray}
I_k(P) &=& P^{-(m-1)m/2} \nonumber\\
&&\times \,\int^{P^{m-1}}_{-P^{m-1}} \int^{P^{m-2}}_{-P^{m-2}}\!\cdots \int^{1}_{-1}
     \Big|\bE\,\exp\{2\pi if (S)\}{\bf 1}_{\{-P<S\le P\}}\Big|^{2k} \,d\alpha_m\cdots d\alpha_1,\nonumber\\
\end{eqnarray}
where ${\bf 1}_{A}$ denote the indicator function of~$A$.
Then for all
  $k=m\tau$ there exists a constant
   $D_{\tau}=D(k, \tau),$ depending on $m$ and $\tau,$ such that one has
\begin{equation}
I_k(P)\le D_{\tau} P^{2k-\triangle(\tau)}
     \Big(\sup_I {\bP}(S\in I)\Big)^{2k},
\end{equation}
where
$$
\triangle(\tau)=0.5 m(m+1)(1-(1-m^{-1})^\tau)
$$
and supremum is taken over all intervals
$I=\{x\cl a<x\le a+1\}$ with arbitrary real~$a$.
\end{theorem}

\remno{3} If $S$ takes values  $1,2,\dots,P$ with equal probabilities then
  $\bE\,\exp\{2\pi if (S)\}$ is a periodic function of period 1 with respect to each coefficient    $\alpha_1,\dots,\alpha_k$. Thus,
$
I_k(P)=2^mP^{-2k} J_k(P),
$
and our result coincides with the original Vinogradov estimate with respect to~$P$.

\remno{4} If $S$ has a uniform distribution on $[-P,P]$, one can get a bound more refined than we obtain from Theorem~8 (see \ [8, Lemmas~9 and~10]).

\begin{theorem} Let $m$ and  $b$~be positive integers and
$b\ge (m+1)/2$. Let $P$~be a real number,
  $P\ge 32$, and $S$~be a random variable with uniform distribution on
   $[-P, P]$. Then for
$k=bm$ we have
$$
I_k(P)\le 2^{5m^2+m} P^{-m^2}\frac{1}{1-m/(2b)}.
$$
\end{theorem}

\remno{5} The proof of Theorem 8 is based on the original proof of the Vinogradov mean value theorem
  (see~[21]). It is more suitable for generalizations to arbitrary random variables with non-lattice support than $p$-adic proof due to Linnik, Karatsuba and Wooley   (see, e.g., [22], [24],~[25]). In particular, Theorem 8 is proved by induction with respect to two parameters
     $P$ and $\tau$ simultaneously. Theorem~10 is the first step of induction when
       $\tau=1$ and~$P$ is arbitrary. This step is trivially clear for lattice case in the original
       Vinogradov mean value theorem.

\begin{theorem} Suppose the conditions of Theorem~$8$ are met. Then for all $k,$ $1\le k\le m,$ we have
$$
I_k(P)\le M_k P^{(3k-1)/2}\Big(\sup_I{\bP}(S\in I)\Big)^{2k},
$$
where $M_k$~is some constant depending on
  $k$ and~$m$ only.
\end{theorem}

\section{Characterization and stability properties of finite quadratic forms}\
Now consider the particular case of functionals -- quadratic forms. The case is very important both in theory and in practice. We start with characterization properties of finite quadratic forms.

Let $Z_1 , \ldots , Z_n$~be  i.i.d. standard normal
random variables and
 $a_1, \ldots , a_n$ be real numbers with
$a_1 + \cdots + a_n \neq 0$ и $a_1^2 + \cdots + a_n^2 \neq 0$.
Suppose
that   $X_1, \ldots , X_n$~are i.i.d. random variables such that
$
a_1 Z_1 + \cdots + a_n Z_n \stackrel{\rm d}{=} a_1 X_1 + \cdots + a_n X_n,
$
where   $\stackrel{\rm d}{=}$ denotes the equality in distribution. Then, by $Cram\acute{e}r's$ decomposition
theorem for the normal law (see\ [27, Theorem~3.1.4]) the random variables $X_i$ are standard normal as well.

Lukacs and Laha (see \ [28, Theorem~9.1.1]) considered а more general problem.
Namely, let   $X_1, \ldots , X_n$~be i.i.d. random variables such that their linear
combination
$L =  a_1 X_1 + \cdots + a_n X_n$ has analytic characteristic function and
$
a_1^s + \cdots + a_n^s \neq 0$ \mbox{for all} $s=1,2,\ldots\,.
$
Then the distribution of $X_1$ is uniquely determined by that of $L$.

The aim of this Section is to obtain а similar characterization property for
quadratic forms in i.i.d. random variables   $Z_1,\ldots,Z_n$. Furthermore, we state а
stability property of such quadratic forms.

Consider а symmetric matrix $A=(a_{ij})_{i,j=1}^n$.  Let
$$
Q (x_1, \ldots , x_n) = \sum^n_{i, j=1} a_{ij} x_i x_j
$$
be а quadratic form in variables $x_1, \ldots , x_n$.
Assume that~$Q$ is non-degenerate
in the sense that   $A$ is not а zero matrix. Suppose  $Z_1,\ldots,Z_n$~are i.i.d. random
variables with а symmetric distribution~$F$.

We say that а pair $(Q,F)$ has а {\it characterization property} (CP), iff for а
sequence of i.i.d. symmetric random variables $X_1, \ldots , X_n$ the equa1ity
\begin{equation}\label {1}
Q (Z_1, \ldots , Z_n) \stackrel{\rm d}{=}   Q (X_1, \ldots , X_n)
\end{equation}
imp1ies
$
Z_1 \stackrel{\rm d}{=} X_1.
$

\remno{6} We require in the definition of CP that the random variables
$X_1,\ldots,X_n$ are symmetric. Otherwise the problem does not have solution even
in the case $n=1$ and $Q(x_1)=x_1^2$. Equation~(\ref{1})
ho1ds for $X_1 = Z_1$, as well as
for $X_1 = |Z_1|.$

\remno{7} With а symmetric distribution $F$ an answer is trivia1 in the
one dimensiona1 case, i.e. any pair $(Q,F)$ has~CP. Therefore we assume that $n\ge 2$ everywhere be1ow.

We are interested to find the sufficient conditions in order the pair   $(Q,F)$ has CP.
The solution of the problem depends a1so on the coefficients of the
matrix~$A$, where the following possibilities exist:

1.~$a_{ii} = 0$ for all $i = 1,\ldots , n$.

2.~$a_{ii} \neq 0$ for some $i = 1, \ldots , n$.

$\phantom{\rm1.~}\vrule width0pt height13pt depth6pt$2.1. $a^{2k+1}_{11} + a^{2k+1}_{22} + \cdots + a^{2k+1}_{nn} \neq 0$ for all $k = 0,1,2,\ldots\,.$

$\phantom{\rm1.~}\vrule width0pt height11pt depth5pt$2.2.~$a_{11} + a_{22} + \cdots + a_{nn} = 0$.

$\phantom{\rm1.~2.2.~}\vrule width0pt height11pt depth5pt$2.2.1.~$a_{ij} = 0$ for all $i \neq j$.

$\phantom{\rm1.~2.2.~}\vrule width0pt height11pt depth5pt$2.2.2.~$a_{ij} \neq 0$ for some $i \neq j$.

$\phantom{\rm1.~}\vrule width0pt height11pt depth5pt$2.3.~$a^{2k+1}_{11} + a^{2k+1}_{22} + \cdots + a^{2k+1}_{nn} = 0$ for some $k = 1,2, \ldots\,.$

\noindent So far the answers are known for the cases~1, 2.1 and~2.2.1 only.

Define now а class   $\ccF$ of probability distributions so that   $F \in\ccF$ iff the
following two conditions are satisfied:

(C1) $F$ has moments $\alpha_k = \int^{\infty}_{-\infty} x^k\,dF(x)$ of all orders $k$;

(C2) $F$ is uniquely specified Ьу   $\alpha_1 , \alpha_2 , \ldots\,.$

\noindent The following examples demonstrate when probability distribution $F$ belongs to~$\ccF$.

\exampno{1}\ If $F$ has an analytic characteristic function, then $F \in {\ccF}$.

Recall (see, e.g., [29, \S\,7.2]), that  а characteristic function is analytic iff

(i) condition (C1) is met and

(ii) $\limsup_{n \to \infty} \alpha_{2n}^{1/(2n)}/(2n) < \infty $.

\noindent The latter condition leads to (C2) (see, e.g., [28, Ch.~9]).

We say that а probability distribution   $F$ satisfies   ${\it Cram\grave{e}r's\,\, condition}$ CC, iff
$$
\int_{-\,\infty}^{\infty} \exp \{ h |x| \} \,dF (x) < \infty \quad
   \mbox{for some} \ h > 0.
$$

Let $F$ satisfies CC, then  $F \in {\ccF}$.

It follows from the fact that   $F$ satisfies  CC iff its characteristic function is
analytic (see, e.g., [29, \S\,7.2]).

\exampno{2}\ If the moments   $\{ \alpha_k \}$ of $F$ satisfy the Carleman condition, i.e.
\begin{equation}\label{c1}
\sum^{\infty}_{n=1} \alpha^{-1/(2n)}_{2n} = \infty,
\end{equation}
then $F \in {\ccF}$.

In fact, the condition   (\ref{c1})
yields the uniqueness of the moment problem for
$\ccF$ (see, e.g., [30, Theorem~1.10]).

Note that the Carleman condition is weaker than  CC. Other examples of probability
distributions belonging to   ${\ccF}$, as well as detailed discussion concerning the
moment problem and other related topics see, e.g.,  [31], [32, Sec.~VII.3] and [33, Sec.~8.12 and~11]).

The following theorem is true (see\ [10, Theorem~3.2.1]).

\begin{theorem}
Let  $F \in {\ccF}$ and the matrix $A$ be such that  $a_{ii} = 0$ for all
$i = 1,\ldots,n$. Then $(Q,F)$ has CP.
\end{theorem}

\exampno{3}\ Let $Z_1 , Z_2 , Z_3$~be i.i.d. standard normal random variables and  $X_1 ,X_2 , X_3$ be i.i.d. symmetric random variables such that
$
Z_1 Z_2 - Z_2 Z_3 \stackrel{\rm d}{=} X_1 X_2 - X_2 X_3.
$
then by Theorem~12 the random variables   $X_1 , X_2 , X_3$~are standard normal.

The following theorem is true (see \ [10, Theorem   3.2.2]).

\begin{theorem}
Let   $F \in {\ccF}$ and the matrix $A$ be such that
$a^{2k+1}_{11} + a^{2k+1}_{22} + \cdots +
a^{2k+1}_{nn} \neq 0$ for all $k = 0,1, 2, \ldots\,.$ Then  $(Q,F)$ has CP.
\end{theorem}

\exampno{4}\ Let $Z_1 , Z_2$~be i.i.d. random variables with distribution~$F$ and
density function
\begin{equation}\label{c2}
p(x) = \frac{1}{4}\, \exp\{-|x|^{1/2}\}, \qquad x \in (-\infty,\infty).
\end{equation}
Then $F \in {\ccF}$ (see\ [33, Ch.~11]).

Let $X_1 , X_2$~be i.i.d. symmetric random variables such that
$
2 Z^2_1+4Z_1 Z_2 - Z^2_2 \stackrel{\rm d}{=}
2 X^2_1+4X_1 X_2 - X^2_2.
$
Then by Theorem 2 the random variables $X_1$ and $X_2$ have the density
function defined in~(\ref{c2}) as well.

The following theorem is true (see\ [10, Theorem~3.2.3]).

\begin{theorem}
   Let $a_{ii} \neq 0$ for some $i\in \{1,\ldots, n\},$ but $a_{11} + \cdots + a_{nn}= 0$ and $a_{ij} = 0$ for all $j \neq i$. Then for any $F$ the pair  $(Q,F)$ does not have~CP.
\end{theorem}

\exampno{5}\ Let $Z$ be а random variable with symmetric distribution   $F$ independent of the random variable $\zeta$ with  ${\bf P}(\zeta =1) ={\bf P}(\zeta =- 1) = 1/2$, and let $c > 0$~be а real constant. Put
$
X = \zeta (Z^2+c)^{1/2}.
$
Suppose now that both
$Z, Z_1, Z_2 , \ldots ,Z_n$ and $X,X_1,X_2,\ldots, X_n$
are
i.i.d. too. Under the conditions of Theorem~13, varying the constant~$c$, we find
а family of symmetric distributions of  $X_1$, such that~(\ref{1}) holds. In particular,
if
$
Z^2_1 - Z^2_2 \stackrel{\rm d}{=} X^2_1 - X^2_2,
$
then the distributions of $X_1$ and $Z_1$ may differ.

Example 5 proves Theorem~13. The proofs of Theorems~11  and~12 are based on the following facts:

a) if $F \in {\ccF}$, then   $X_1$ has moments  ${\bf E}\, X^k_1$ of all orders~$k$;

b) under the given conditions, we have
$
{\bf E} \, X^k_1 = {\bf E} \, Z^k_1$ \mbox{for all} $k =
1,2, \ldots \,.
$

The following stability theorem is true as well (see \ [10, Theorem~3.2.4]).

\begin{theorem} Suppose that the pair $(Q,F)$ has $CP$. Let $X_{N,1},
\ldots, X_{N,n}$ for $N = 1,2,\ldots,$~be а series оf i.i.d. symmetric random variables and
$$
Q (X_{N,1}, \ldots , X_{N,n}) \stackrel{\rrd}{\longrightarrow}  Q (Z_1,
\ldots , Z_n) \quad \mbox{as} \  N \to \infty,
$$
where   $\stackrel{\rm d}{\longrightarrow}$ denotes the convergence in distribution. Then
$$
X_{N,1}  \stackrel{\rrd}{\longrightarrow} Z_1 \quad \mbox{as} \ N \to\infty.
$$
\end{theorem}

We use the tightness of the converging
sequences of quadratic forms while proving Theorem~14.

\section{Distributions of quadratic forms in Gaussian random variables}\
Let ${\sH}$ be a real separable Hilbert space with inner product $(\cdot, \cdot)$ and norm $|x| =
(x,x)^{1/2} $, $x \in {\sH}$. We denote by $Y$ an ${\sH}$-valued
Gaussian random element with zero mean and covariance operator~$V$, i.e., the characteristic
functional of~$Y$ has the form
$$
\bE\,\exp\{i(x,Y)\} = \exp \bigg(-\frac{(Vx,x)}{2}\bigg), \qquad x\in \sH,
$$
where $V$ is defined for all $x, z \in \sH$ by $(Vx,z) = \bE\,(Y, x)(Y,z)$. Assume that the eigenvalues $\sigma^2_j$, $j = 1, 2, \dots,$ of the covariance operator $V$ are ordered so that
$$
\sigma^2_1 = \sigma^2_2 = \cdots = \sigma^2_k > \sigma^2_{k+1}\ge \sigma^2_{k+2}\ge \cdots
$$
with some $k \ge 1$, which is the multiplicity of $\sigma^2_1$.
It is well known that there exists an orthonormal system $\{e_i\}^{\infty}_{i=1}$ in ${\sH}$ such that $Ve_i = \sigma^2_ie_i$, or $Y = \sum^{\infty}_{i=1}Y_i e_i$, where $Y_i = (Y, e_i) $ are independent real Gaussian random variables
with zero mean and variances $\sigma_i^2$. The density function of the random variable $|Y - a|^2$ will be denoted by
$p(u, a)$, $u\ge 0$, $a \in \sH$. Let $f_k(u, a)$ be a density function of $\sum^{k}_{i=1}(Y_i - a_i)^2$.
We recall (see, e.g., [34, vol.~6, p.~421]) that the exact expression for the density
  $f_k(u, a)$ is known.

Set $R = \sum^{\infty}_{j=k+1}(Y_j - a_j)^2$. Note that $\bE\, R = \sum^{\infty}_{j=k+1}(\sigma^2_j + a_j^2).$

Let $u_0 = 2k(1-\sigma^2_{k+1}\sigma^{-2}_1)^{-2} \,\bE\, R$, $\overline a_i=(a_1,\ldots,a_i)$.

The following theorem is true (see \ [11, Theorem 1]).

\begin{theorem} \label{th4.1}
Suppose $k \ge 3$.

{\rm a)} For all $u\ge 0$ и $a\in \sH$ the upper bound
\begin{equation}
p(u, a) \le f_k(u, a) \,\bE\,\exp \bigg\{\frac{R}{2\sigma^2_1}\bigg\},
\end{equation}
holds, where
\begin{equation}
 \bE\, \exp \bigg\{\frac{R}{2\sigma^2_1}\bigg\} = \prod^{\infty}_{j=k+1}\bigg(1 - \frac{\sigma^2_j}{\sigma^2_1}\bigg)^{-1/2}\exp\bigg\{\frac{a^2_j}{2(\sigma^2_1 - \sigma^2_j)}\bigg\}.
\end{equation}

{\rm b)} Let one of the following conditions~{\rm i)--ii) be met:}

{\rm i)} $k=3,$ $u \ge u^* = 4.9 u_0 + 16.94|\overline{a}_3|^2\sigma^{-4}_1u^2_0$\kbt

{\rm ii)} $k\ge 4,$ $u \ge u^{**} = 5.625(k-1)^2(k-3)^{-1}u_0 + 16|\overline{a}_3|^2\sigma^{-4}_1u^2_0(k-3)^{-2}.$
Then   the lower bound
\begin{equation}
p(u, a) \ge 0.125 f_k(u, a) \,\bE\,\exp \bigg\{\frac{R}{2\sigma^2_1}\bigg\}.
\end{equation}
holds.
\end{theorem}

It follows from (17) and (18) that the properties of  $p(u, a)$ with respect to $u$
and $\overline{a}_k$ could be described via properties of $f_k(u, a)$.
For example, it was proved in~[35] (see \ ~(15) in~[35])
the following inequality: for all
  $k\ge 1$, $u > 0$ and arbitrary $a\in \sH$ one has
\begin{equation}
f_k(u, a) \le \bigg(2\sigma^2_1\Gamma\bigg(\frac{k}{2}\bigg)\bigg)^{-1}
\bigg(\frac{u}{2\sigma_1^2}\bigg)^{k/2-1}
\exp\bigg\{-\frac{(u^{1/2}-|\overline{a}_k|)^2}{2\sigma^2_1}\bigg\}.\nonumber
\end{equation}
In Lemma 2 in [36] another bound was obtained: if $k\ge 1$ and $|\overline{a}_k| > 0$, then for all $u > 0$ one has
\begin{equation}
f_k(u, a) \le c_1(k)\sigma^{-1}_1 u^{(k-3)/4}
|\overline{a}_k|^{-(k-1)/2}
\exp\bigg\{-\frac{(u^{1/2}-|\overline{a}_k|)^2}{2\sigma^2_1}\bigg\},\nonumber
\end{equation}
where $c_1(k) = \pi^{-1/2} + ((k-1)/2)^{(k-1)/2}/\Gamma(k/2)$ and $\Gamma(\cdot)$~is the Gamma function.

Furthermore, since
$$
{\bP}(|Y-a| > r) = \int_{r^2}^{\infty}p(u,a)\,du,
$$
the two-sided bounds for   $p(u,a)$ imply two-sided bounds for the probability   ${\bP}(|Y-a| > r)$.
In particular,
the following theorem is true (see \ [11, Corollary~4,i)]).

\begin{theorem} \label{col4.1}
Let $k\ge 4$ and $r,a$ be such that $|\overline{a}_k| > 0$ and $r > \sigma_1^2/|\overline{a}_k| + 2|\overline{a}_k| + (u^{**})^{1/2}$. Then there exist constants $c_9(k), c_{10}(k),$ depending on $k,$ only such that one has
\begin{eqnarray*}
c_9(k)\,\bE\,\exp \bigg\{\frac{R}{2\sigma^2_1}\bigg\}& \le& {\bP}(|Y-a| > r) \exp \bigg\{\frac{(r - |\overline{a}_k|)^2 }{2\sigma^2_1}\bigg\}\sigma_1r^{(k-3)/2}|\overline{a}_k|^{(k-1)/2}\\
& \le& c_{10}(k)\,\bE\,\exp \bigg\{\frac{R}{2\sigma^2_1}\big\}.
\end{eqnarray*}
\end{theorem}

Another type of bounds for ${\bP}(|Y-a| > r)$ see \ in~[11] and~[37], and in the papers from bibliographies in~[11] and~[37].


\begin{thebibliography}{99}

\bibitem{art2}
{\it G\"{o}tze~F.} Asymptotic expansions for bivariate von Mises functionals.~---
Z.\ Wahr\-scheinlichkeitstheor.\ verw.\ Geb., 1979, v.~50, No 3, p.~333-–355.

\bibitem{art20}
{\it G\"{o}tze~F.} Edgeworth expansions in functional limit theorems.~--- Ann.\ Probab., 1989, v.~17, No 4, p.~1602--1634.

\bibitem{art3}
{\it Weyl~H.} \"{U}ber die Gleichverteilung von Zahlen mod.\ Eins.~--- Math.\ Ann., 1916, v.~77, No 3, p.~313--352.

\bibitem{art1}
{\it Yurinskii~V.\,V.} An estimate of error of the normal approximation of the probability of a hit
into a ball.~--- Dokl. Akad. Nauk SSSR, 1981, т.~258, с.~577--578 (in Russian).

\bibitem{art100}
{\it Yurinskii~V.\,V.} On an error of the normal approximation,.~--- Sibirsk.\ Mat.\
Zh., 1983, т.~24, No 6, c.~188--199 (in Russian).

\bibitem{art101}
{\it G\"{o}tze~F.\kb Prokhorov~V.\,V.\kb Ulyanov~V.\,V.} Bounds for characteristic functions of polynomials
in asymptotically normal random variables.~--- Russian Math.\ Surveys, 1996, т.~51, No 2, c.~181--204.

\bibitem{art102}
{\it G\"{o}tze~F.\kb Prokhorov~V.\,V.\kb Ulyanov~V.\,V.} On smooth behavior of probability distributions under polynomial mappings.~--- Theory Probab.\ Appl., 1998, Vol.~42, No.~1, p.~28--38.

\bibitem{art1021}
{\it G\"{o}tze~F.\kb Prohorov~V.\,V.\kb Ulyanov~V.\,V.} A stochastic analogue of the Vinogradov mean value theorem.~--- Probability Theory and Mathematical Statistics (Tokyo, 1995). River Edge: World Scientific, 1996, p.~62--82.

\bibitem{art1022}
{\it Prokhorov~Yu.\,V.\kb Christoph~G.\kb Ulyanov~V.\,V.} On Characteristic Properties of Quadratic Forms.~--- Doklady Mathematics, 2000, Vol.~62, No 1, p.~39--41.

\bibitem{art1023}
{\it Christoph~G.\kb Prohorov~Yu.\kb Ulyanov~V.} Characterization and stability problems for finite quadratic forms.~--- Asymptotic Methods in Probability and Statistics with Applications (St.\ Petersburg, 1998). Boston: Birkh\"auser, 2001, p.~39--50.

\bibitem{art1024}
{\it Christoph~G.\kb Prokhorov~Yu.\,V. \kb Ulyanov~V.\,V.}  On distribution of quadratic forms in Gaussian random variables.~--- Theory Probab.\ Appl., 1996, Vol.~40, No.~2, p.~250--260.

\bibitem{art1026}
{\it Prokhorov~Yu.\,V.\kb Ulyanov~V.\,V.} Some approximation problems in statistics and probability.~--- Limit Theorems in Probability, Statistics and Number Theory. Heidelberg: Springer, 2013, p.~235--249. (Springer Proc.\ Math.\ Statist., v.~42.)

\bibitem{art1025}
{\it Prokhorov~Yu.\,V.} Selected Works. Moscow: Torus-Press, 2012, 778~p. (in Russian).

\bibitem{art103}
{\it Rudin~W.} Fourier Analysis on Groups. New York: Interscience Publ., 1962, 285~p.

\bibitem{art104}
{\it Gel'fand~I.\,M.\kb Raikov~D.\,A.\kb Shilov~G.\,E.} Commutative normed rings.~--- Uspekhi Mat.\ Nauk, 1946, Vol.~1, No 2, p.~48--146 (in Russian).

\bibitem{art108}
{\it Zygmund,~A.} Trigonometric Series, vol. I and II combined, Cambridge University Press, London, Cambridge, 1968.
\bibitem{art109}
{\it Wiener~N.\kb Wintner~A.} Fourier--Stieltjes transforms and singular infinite con\-volutions.~--- Amer.\ J.\ Math., 1938, v.~60, No 3, p.~513--522.

\bibitem{art110}
{\it Salem~R.} On singular monotonic functions of the Cantor type.~--- J.\ Math.\ Phys.\ Mass.\ Inst.\ Tech., 1942, v.~21, p.~69--82.

\bibitem{art111}
{\it Schaeffer~A.\,C.} The Fourier--Stieltjes coefficients of a function of bounded variation.~--- Amer.\ J.\ Math., 1939, v.~61, p.~934--940.

\bibitem{art105}
{\it Prachar~K.} Primzahlverteilung. Berlin: Springer-Verlag, 1978, 421~p.

\bibitem{art106}
{\it Vinogradov~I.\,M.} Selected Works, Springer-Verlag, Berlin, 1985, 401~p.

\bibitem{art107}
{\it Arkhipov~G.\,I.\kb Karatsuba~A.\,A.\kb Chubarikov~V.\,N.} The Theory of Multiple Trigonometric
Sums, Nauka, Moscow, 1987, 368~p. (in Russian).

\bibitem{art113}
{\it Stechkin~S.\,B.} On mean values of the modulus of a trigonometric sum.~--- Proc.\ Steklov Inst.\ Math.\ vol.~134, 1975, p.~321--350.

\bibitem{art114}
{\it Wooley~T.\,D.} Vinogradov's mean value theorem via efficient congruencing.~--- Ann.\ Math., 2012, v.~175, No 3, p.~1575--1627.

\bibitem{art115}
{\it Parsell~S.\,T.\kb Prendiville~S.\,M.\kb Wooley~T.\,D.}  Near-optimal mean value estimates for multidimensional Weyl sums.~--- Geom.\ Funct.\ Anal., 2013, v.~23, No 6, p.~1962--2024.

\bibitem{art116}
{\it Linnik~U.\,V.} On Weyl's sums.~--- Матем.\ сб., 1943, т.~12, No 1, p.~28--39.

\bibitem{art117}
{\it Linnik~U.\,V.\kb Ostrovskii~I.\,V.} Decomposition of random variables and vectors. Translations of Mathematical Monographs, Vol. 48, American Mathematical Society, Providence, R. I., 1977, 380~p.

\bibitem{art118}
{\it Lukacs~E.\kb Laha~R.\,G.} Applications of characteristic functions. London: Griffin, 1964, 202~p.

\bibitem{art119}
{\it Lukacs~E.} Characteristic functions, 2nd ed., London: Griffin, 1987, 360~p.

\bibitem{art120}
{\it Shohat~J.\,A.\kb Tamarkin~J.\,D.}  The Problem of Moments. Providence:
Amer.\ Math.\ Soc., 1963, 144~p.

\bibitem{art1201}
{\it Akhizer~N.\,I.} The classical Moment Problem. Oliver\&Boyd, Edinburgh, 1965.


\bibitem{art121}
{\it Feller~W.} An Introduction to Probability Theory and its Applications II. New York: Wiley, 1971.

\bibitem{art122}
{\it Stoyanov~J.\,M.} Counterexamples in Probability, New York: Wiley, 1987.

\bibitem{art10}
Encyclopedia of Mathematics (Transl. of the Soviet "Mathematical Encyclopedia"), Vol. 6,
Kluwer Academic Publishers, Norwell, MA, 1990.

\bibitem{art123}
{\it Ulyanov~V.\,V.} Nonuniform estimates of the density of the squared norm of a Gaussian vector
in Hilbert space.~--- New Trends in Probability and Statistics, v.~1. Ed.\ by V.\,Sazonov and T.\,Shervashidze. Utrecht/Vilnius: VSP/Mokslas, 1991, p.~255--263.

\bibitem{art4}
{\it Ulyanov~V.\,V.} On Gaussian measure of balls in $H$.~--- Frontiers in Pure and Applied Probability. II: Proceedings of the Fourth Russian--Finnish Symposium on Probability
Theory and Mathematical Statistics. Ed.\ by A.\,N.\,Shiryaev et al. Moscow:
TVP Science Publishers, 1996, p.~195--206.

\bibitem{art41}
{\it Rozovsky~L.\,V.} On Gaussian measure of balls in a Hilbert space.~--- Theory Probab.\ Appl., 2009, v.~53, No 2, pp.~357--364.

\end{thebibliography}
\end{document}